\documentclass[11pt]{article}
\usepackage {graphics}
\usepackage {pstricks,pst-node}
\usepackage{amsmath,amssymb,amsfonts,oldlfont,latexsym,dsfont,amsthm}
\usepackage[T1]{fontenc}
\newtheorem{tw}{Theorem}

\newtheorem{wn}[tw]{Corolary}

\newtheorem{hip}[tw]{Conjecture}

\onecolumn
\begin{document}

\newcommand{\trou}{\vspace{5 mm} \noindent}
\newcommand{\Trou}{\vspace{7 mm} \noindent}
\newcommand{\tour}{\vspace{4 mm}}
\newcommand{\ptrou}{\vspace{4 mm} \noindent}
\newcommand{\cbdo}{\hfill \rule{.1in}{.1in}}
\newcommand{\prf}{\noindent{\bf Proof.}\ }
\newcommand{\cex}{\noindent{\bf Counterexample.\ }}

\newcommand{\n}{I\!\!N}
\newcommand{\z}{{\sf Z}\!\!{\sf Z}}
\newcommand{\q}{Q\!\!\!\!I}
\newcommand {\rn}[1]{{I\!\!R}^{#1}}
\newcommand{\ce}{C\!\!\!\!I}
\newcommand {\cen}[1]{{C\!\!\!\!I}\,^{#1}}
\newcommand{\pusty}{\mbox{\sc\O}}

\newcommand{\s}{\color{red}}


\title{A note on edge {colorings}\\ distinguishing
 all triangles in a graph}
\author{Monika Pil\'sniak and Mariusz Wo\'zniak}
\date{\small Department of Discrete Mathematics\\
AGH University, al. Mickiewicza 30, 30-059 Krak\'ow, Poland}

\maketitle
\begin{abstract}

We consider edge colorings of a graph in such a way that each two different triangles have distinct colorings. It is {an extension} of the well-known idea of distinguishing all maximal stars in a graph. It was introduced in literature in 1985 and studied by many authors in various variants, but always for stars. We estimate new invariants {regarding triangles} for proper and general colorings.
\end{abstract}

\section{Introduction}
\par

Let $G=(V, E)$ be a finite graph with a given edge coloring. Consider a maximal star $K_{1, {\rm{deg}} (v)}$ at every vertex $v$.  Two stars $S_1$ and $S_2$ with $|V(S_1)| \geq |V(S_2)|$ are {\it distinct} if either there exists an edge of $S_1$ with color $\alpha$ and there is no edge of $S_2$ with color $\alpha$ or there exists a color $\beta$ such that the number of edges of $S_1$ with $\beta$ is different from the number of edges with $\beta$  in $S_2$ (the coloring does need not be proper). An edge coloring is called a {\it coloring distinguishing stars} (in the literature called also vertex distinguishing) if for any two different vertices $v$ and $w$ the stars $K_{1, {\rm{deg}} (v)}$ and $K_{1, {\rm{deg}} (w)}$ are distinct.

In 1996, Burris and Shelp \cite{BuSc97}, and independently \v{C}ern\'{y}, Hor\v{n}\'{a}k and Soták \cite{CHS96} posed a question: {\it How many colors are needed to distinguish stars in a proper coloring of a graph?}. It was conjectured by Burris and Schelp and proved in \cite{MW1}, that if $G$ is a graph of order $n$, without isolated edges and with at most one isolated vertex, then $n + 1$ colors enough to distinguish all maximal stars in $G$. 
This topic was spanned by hundreds of papers and was in the main spectrum of interest of {many mathematicians}. It was based on a similar {concept} introduced ten years earlier by  Harary and Plantholt \cite{HaPl85} for general vertex colorings (not necessary proper) and later by Aigner and Triesch \cite{AiTr90a}.

\trou

This note is the first paper where another graph is investigated instead of a star.
More precisely, we will consider triangles $K_3$. Let $f:E\mapsto \{1,2,\ldots, k\}$ be an edge coloring of a graph $G=(V,E)$. Then, each triangle $T\subset G$ defines a set (or a multiset) of colors
$F(T):=\{f(e):e\in E(T)\}$ called a \emph{palette}.

We say that a coloring $f$ \emph{distinguishes triangles} if for every two different triangles $T_1$, $T_2$
we have $F(T_1)\neq F(T_2)$.

We will only consider the case $G=K_n$, $n\ge 3$.
However, observe  that the obtained results also give estimates of the number of colors needed for any graphs.
Indeed, if a coloring distinguishes triangles in a complete graph, then this property is preserved when any number of edges is removed.

In Section \ref{s_proper} we investigate a proper coloring $f$ and determine  the smallest number $\tau'(n)$ of colors such that the coloring $f$ distinguishes triangles in $K_n$.

In Section \ref{s_general} we investigate a general (not necessary proper) coloring $f$. Then the palettes are multisets, and by $\tau_m(n)$ we denote the smallest number $k$ of colors in an edge coloring $f$ distinguishing triangles in $K_n$. We estimate $\tau_m(n)$ and we discuss a modification $\tau_s(n)$, when we consider sets of colors of edges instead of multisets, despite the fact that colorings need not be proper.

\section{Proper coloring}\label{s_proper}

\begin{tw}
\[\tau'(n)=
\left\{
\begin{array}{ll}
n & \mbox{if $n$ is odd} \\
n+1 & \mbox{if $n$ is even. }
\end{array}
\right.
\]
\end{tw}
\noindent{{\bf Proof.}} Note first that since $f$ is proper, the pallets {of triangles}
$F(T)$ are \emph{rainbow}, \emph{i.e.} all colors are distinct. Hence,
if we use $k$ colors, we can get at most
${k\choose 3}$  different palettes. On the other hand, $K_n$ contains ${n\choose 3}$ triangles.
{So $k\ge n$}.

Suppose {$n$ is odd}. Let us denote the vertices of the graph by
$x_1, x_2, \ldots, x_n$ and assign color $i+j \pmod n$ to the edge $x_ix_j$.
It is easy to see that the coloring defined in this way is proper. Now,
it is enough to show that every possible triplet of colors is present as a palette
of some triangle.
Let $a,b,c$ be three  different colors. We are looking for a solution to the system {of equations}
$$
\begin{array}{lll}
p + q & = & a \\
q + r & = & b \\
r + p & = & c. \\
\end{array}
$$
where the addition is taken  $\mod n$. Such a system always has exactly one solution such that the numbers $p, q, r$ are pairwise different, which means that the vertices $x_p, x_q, x_r$ define
`a real' triangle. To show this, it may be the easiest way to give \emph{explicit} formulas for $p, q, r$.

We have
$$\begin{array}{lll}
p & = & (\frac{n+1}{2})(a-b+c) \\
q & = & (\frac{n+1}{2})(a+b-c)\\
r & = & (\frac{n+1}{2})(-a+b+c).\\
\end{array}
$$

Now suppose that {$n$ is even}. Let us show that a coloring with $k=n$ colors cannot distinguish all triangles. For then all three different colors would have to be present.
We have exactly $n-1 \choose 2$ triplets containing color `1'.
Let $m_1$ be the number of edges with color `1'. On the other hand we have  $n-2$ triangles in which a color of one edge is fixed. Thus, at most $m_1 (n-2)$ triangles
 can be built
 on  edges colored with `1'. Hence
$$ m_1 (n-2) \ge {n-1 \choose 2} \ge \frac{1}{2} (n-1)(n-2), $$
which gives $m_1\ge \frac{1}{2} (n-1)$. Since $n$ is even we have to have  $m_1=\frac{1}{2}n$. However, if $m_1=\frac{1}{2}n$,
then on $m_1$ edges
we would have exactly $\frac{1}{2}n(n-2)$ triangles, which is more than there are triplets containing color `1', a contradiction.

However, it is easy to see that there is an appropriate coloring of $K_n$  with $n+1$ colors. It is enough to use the coloring defined in the first part of the proof for $K_{n+1}$. Notice that then $n+1$ is odd. Next, from the graph
$K_{n+1}$ we remove one vertex and we get a coloring $K_n$ with the desired properties. \cbdo

\begin{wn}
If $G$ is connected graph  of order  $n\ge 4$, then $\tau'(G)\leq n$ if $n$ is odd, and $\tau'(G)\leq n+1$ if $n$ is even.
\end{wn}

\section{General coloring}\label{s_general}

If the coloring is general, then the palettes of three colors may be multisets.
Then color palettes do not have to be rainbow, they can be of the form $[\alpha\alpha\beta]$ or $[\alpha\beta\beta]$ or
$[\alpha\alpha\alpha]$. Obviously, by the result for proper colorings, $\tau_m(G)\leq n+1.$  Let us count how many different triples we can get with $k$ colors.
Simple calculations (or using Pascal's triangle property) give

$${k \choose 3} + 2 {k \choose 2} + k = {k+2 \choose 3}.$$

The above formula shows that having $k$ colors at our disposal, we could try to obtain a coloring that distinguishes all triangles of the complete graph $K_n.$ for $n=k+2$. This is the case, for example, for $n=3$,
where of course one color is enough.
The next theorem shows that this is the only exception.

\begin{tw}
For $n\ge 4$, we have $\tau_m(n)\ge n-1.$
\end{tw}

\noindent{{\bf Proof.}} Suppose there is a coloring $f$ that uses $k=n-2$ colors.
As we have shown above, each 3-element multiset of colors has its realization as a palette of exactly one
triangle in $K_n$. We will count how many edges colored with $\alpha$ are needed to complete all pallets
as multisets of the form $[\alpha,\beta,\beta]$.

Note that an edge colored with $\alpha$ can {belong only to} one triangle with a palette
$[\alpha,\beta,\beta]$. For, suppose that an edge $e=x_1x_2$ colored with $\alpha$ belongs to two triangles
$T$ and $T'$, where $V(T)=\{x_1x_2y\}$ and $V(T')=\{x_1x_2z\}$. Suppose, moreover, that
$f(yx_1)=f(yx_2)=\beta$, and $f(zx_1)=f(zx_2)=\gamma$. Then two triangles of $K_n$ with vertices
$x_1,y,z$ and $x_2,y,z$ have the same color palette.

The number of multisets of the form $[\alpha,\beta,\beta]$ with a fixed color $\alpha$  for $\beta\neq \alpha$
is $k-1$. In addition, to implement a multiset $[\alpha,\alpha,\alpha]$ we need three edges. Thus, to realize all multisets of the form $[\alpha,\beta,\beta]$ we need $k+2$  edges of color $\alpha$. In total, to realize all multisets of the form $[\alpha,\beta,\beta]$, where $\alpha$ is one of $k$ colors,
we need at least $k(k+2)$ edges. This number must not be greater than the number of all edges. Hence we have

$$k(k+2)\le {n \choose 2} =  {k+2 \choose 2}.$$

The above inequality also holds only for $k=1$, \emph{i.e.} for $n=3$, a contradiction. \cbdo

\trou
Perhaps the following conjecture holds.

\begin{hip}\label{hip_multi}
For $n\ge 4$ we have $\tau_m(n)= n-1.$
\end{hip}

\Trou

Another possibility is to distinguish triangles by sets in a general coloring.
There are fewer 3-element sets than multisets,
but more than 3-element rainbow sets. Hence, we have an obvious inequality
$$\tau_m(n)\le \tau_s(n)\le \tau'(n).$$

The following conjecture would imply Conjecture~\ref{hip_multi}.

\begin{hip}\label{hip_zbiory}
For $n\ge 4$ we have $\tau_s(n)= n-1.$
\end{hip}

\begin{figure}[htb]
\psset{unit=1cm}
\psset{radius=0.2}

\begin{pspicture}(12,3.5)

\put(1,0)
{\begin{pspicture}(0,0)
\dotnode(0,0){A}
\dotnode(3,0){B}
\dotnode(1.5,2.5){C}
\dotnode(1.5,1){X}
\ncline{B}{A}
\ncline{B}{C}
\ncline{C}{A}
\ncline{X}{A}
\ncline{B}{X}
\ncline{C}{X}

\rput(1.5,1.75){\circlenode[linestyle=none,fillstyle=solid,fillcolor=green]{C1}{\mbox{\tiny 1}}}
\rput(0.75,0.5){\circlenode[linestyle=none,fillstyle=solid,fillcolor=green]{C1}{\mbox{\tiny 2}}}
\rput(2.25,0.5){\circlenode[linestyle=none,fillstyle=solid,fillcolor=green]{C1}{\mbox{\tiny 3}}}
\rput(1.5,0){\circlenode[linestyle=none,fillstyle=solid,fillcolor=green]{C1}{\mbox{\tiny 1}}}

\rput(0.75,1.25){\circlenode[linestyle=none,fillstyle=solid,fillcolor=green]{C1}{\mbox{\tiny 1}}}
\rput(2.25,1.25){\circlenode[linestyle=none,fillstyle=solid,fillcolor=green]{C1}{\mbox{\tiny 1}}}

\end{pspicture}
}

\put(9,0)
{\begin{pspicture}(0,0)
\dotnode(-1,0){A}
\dotnode(1,0){B}
\dotnode(1.5,1.5){C}
\dotnode(0,2.5){D}
\dotnode(-1.5,1.5){E}
\dotnode(0,3.5){X}
\ncline{B}{A}
\ncline{B}{C}
\ncline{C}{D}
\ncline{E}{D}
\ncline{E}{A}

\ncline{C}{A}
\ncline{D}{B}
\ncline{C}{E}
\ncline{D}{A}
\ncline{E}{B}

\ncline{C}{X}
\ncline{D}{X}
\ncline{E}{X}

\pscurve(1,0)(2,1)(0,3.5)
\pscurve(-1,0)(-2,1)(0,3.5)
\rput(-2,1){\circlenode[linestyle=none,fillstyle=solid,fillcolor=green]{C1}{\mbox{\tiny 4}}}
\rput(2,1){\circlenode[linestyle=none,fillstyle=solid,fillcolor=green]{C1}{\mbox{\tiny 3}}}

\rput(0,0){\circlenode[linestyle=none,fillstyle=solid,fillcolor=green]{C1}{\mbox{\tiny 4}}}
\rput(0,1.5){\circlenode[linestyle=none,fillstyle=solid,fillcolor=green]{C1}{\mbox{\tiny 3}}}

\rput(-0.5,0.4){\circlenode[linestyle=none,fillstyle=solid,fillcolor=green]{C1}{\mbox{\tiny 5}}}
\rput(0.5,0.4){\circlenode[linestyle=none,fillstyle=solid,fillcolor=green]{C1}{\mbox{\tiny 1}}}

\rput(-0.5,1.2){\circlenode[linestyle=none,fillstyle=solid,fillcolor=green]{C1}{\mbox{\tiny 2}}}
\rput(0.5,1.2){\circlenode[linestyle=none,fillstyle=solid,fillcolor=green]{C1}{\mbox{\tiny 4}}}

\rput(-1.25,0.75){\circlenode[linestyle=none,fillstyle=solid,fillcolor=green]{C1}{\mbox{\tiny 5}}}
\rput(1.25,0.75){\circlenode[linestyle=none,fillstyle=solid,fillcolor=green]{C1}{\mbox{\tiny 3}}}

\rput(-0.75,2){\circlenode[linestyle=none,fillstyle=solid,fillcolor=green]{C1}{\mbox{\tiny 1}}}
\rput(0.75,2){\circlenode[linestyle=none,fillstyle=solid,fillcolor=green]{C1}{\mbox{\tiny 2}}}

\rput(0,3){\circlenode[linestyle=none,fillstyle=solid,fillcolor=green]{C1}{\mbox{\tiny 1}}}

\rput(-0.75,2.5){\circlenode[linestyle=none,fillstyle=solid,fillcolor=green]{C1}{\mbox{\tiny 5}}}
\rput(0.75,2.5){\circlenode[linestyle=none,fillstyle=solid,fillcolor=green]{C1}{\mbox{\tiny 2}}}

\end{pspicture}
}

\end{pspicture}
\caption {Edge colorings showing that
$\tau_m(4)= \tau_s(4) = 3$, and $\tau_m(6)=\tau_s(6)=5$} \label{rys} 
\end{figure}
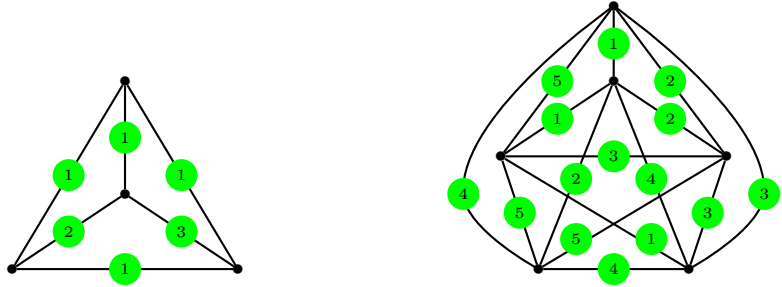

Some colorings supporting the above conjectures are shown in Figure~\ref{rys}.

\trou

\noindent{\bf Acknowledgment.} The authors would like to thank our
colleague Jakub Przyby{\l}o for many helpful comments.

\trou

\end{document}